\documentclass[12pt]{article}
\begin{document}
\newcommand{\qed}{\ \
\mbox{\rule{8pt}{8pt}}\vspace{0.3cm}\newline}
\newcommand{\ia}{{\bf I}_{A_i}}
\newcommand{\ba}{\widetilde{{\bf KC}}_i}
\newcommand{\bba}{{\bf KC}_i}
\newcommand{\ra}{\longrightarrow}
\newcommand{\pe}{{\cal P}}
\newcommand{\der}{{\cal DP}_d}
\newcommand{\ot}{\otimes}
\newcommand{\rec}{\raisebox{-1ex}{\ $\stackrel{\textstyle{\stackrel{\textstyle{\longleftarrow}}{\longrightarrow}}}{\longleftarrow}$\ }}
\title{Derived Kan Extension for strict polynomial functors}
\author{Marcin Cha\l upnik \\
\normalsize{Institute of Mathematics, University of Warsaw,}\\
\normalsize{ul.~Banacha 2, 02--097 Warsaw, Poland;}\\
\normalsize{Institute of  Mathematics PAN,}\\
\normalsize{ul.~\'Sniadeckich 8, 00--956 Warsaw, Poland.}\\
\normalsize{e--mail: {\tt mchal@mimuw.edu.pl}}}
\date{\mbox{}}
\newtheorem{prop}{Proposition}[section]
\newtheorem{cor}[prop]{Corollary}
\newtheorem{theo}[prop]{Theorem}
\newtheorem{lem}[prop]{Lemma}
\newtheorem{defi}[prop]{Definition}
\newcommand{\ka}{{\mbox {\bf k}}}
\maketitle
\begin{abstract}
We investigate fundamental  properties of adjoint functors to the precomposition functor in the category of strict polynomial functors.
In particular, we obtain the proof of Touz\'e's Collapsing Conjecture on Ext--groups between twisted functors  and refinement of a general program of computing
Ext--groups between strict polynomial functors important in representation theory.\\\mbox{}\vspace{0.1cm}
\end{abstract}
\section{Introduction}
The category ${\cal P}_d$ of homogeneous strict polynomial functors of degree $d$ over a finite field \ka\ has proved to be a valuable tool
in investigating homological problems concerning $\ka[GL_n(\ka)]$--modules. For example the proof of finite generation of the rational cohomology
ring of a finite group scheme \cite{FS} heavily depended on construction and properties  of certain ``universal classes'' in Ext--groups in ${\cal P}_d$.
The treatment of cohomology rings for reductive groups in \cite{TvdK} depends on similar universal classes constructed by Touz\'e \cite{T1,T3}.\newline
The most striking feature of ${\cal P}_d$, which makes it so useful in applications is relative ease (when compared to related module categories) with which one can compute the Ext--groups there. These computations (see e.g. \cite{FS, FFSS, C2, C3, C4}) have focused mainly on the Ext--groups of the form
$\mbox{Ext}^*_{{\cal P}_{dp^i}}(F^{(i)},G^{(i)})$ ($F^{(i)}, G^{(i)}$ mean the i--th Frobenius  twists of $F,G\in {\cal P}_d$), which are closely related to the Ext--groups in the category   $\ka[GL_n(\ka)]$--mod by \cite{FFSS, Be}.
The main part of all these works was comparing $\mbox{Ext}^*_{{\cal P}_{dp^i}}(F^{(i)},G^{(i)})$  to  Ext--groups between some
untwisted functors which were easier to handle. This was usually achieved by considering hyper--Ext groups for complexes connecting twisted and untwisted functors and required massive calculations of spectral sequences.
However, it was recently observed  by Touz\'e \cite{T2} that this seemingly complicated process of  ``untwisting the Ext--groups'' fits into a nice  general scheme. He constructed for any $F,G\in{\cal P}_d$ a spectral sequence converging to  $\mbox{Ext}^*_{{\cal P}_{dp^i}}(F^{(i)},G^{(i)})$, with $E_2$--page consisting of
 Ext--groups between untwisted functors. He also proved that this spectral sequence collapses at $E_2$ in some special cases and formulated the ``Collapsing Conjecture''
stating that it collapses at $E_2$ for all $F,G\in{\cal P}_d$. The main result of the present paper is a proof of a stronger version of this conjecture. Namely we show \newline\mbox{}\vspace{0.5cm}\newline
{\bf Corollary 3.7} {\em For any $F,G\in {\cal P}_d$ and $i\geq 0$, there is a natural in $F,G$ isomorphism of graded spaces
\[\mbox{Ext}^*_{{\cal P}_{dp^i}}(F^{(i)},G^{(i)})\simeq \mbox{Ext}^*_{{\cal P}_d}(F,G\circ(-\otimes A_i))\]
where $A_i$ is a graded space which is one--dimensional in nonnegative even degrees strictly smaller than $2p^i$ ($p$ is the characteristic of the  ground field)
and trivial elsewhere.}\newline\mbox{}\vspace{0.5cm}\newline
This fact turns out to follow from quite a general phenomenon we call the Derived Kan Extension  and investigate in the present paper. Namely, for $H\in{\cal P}_s$ we construct
the functors  ${\bf K}^r_H, {\bf K}_H^l: {\cal DP}_{ds}\longrightarrow {\cal DP}_{d}$ between derived categories which are respectively the right and left adjoint to the functor of taking precomposition with $H$. The Collapsing Conjecture immediately follows from a computation of the unit of adjunction
 $({\bf C}_H, {\bf K}^r_H)$ for
$H=I^{(i)}$ we provide in Theorem 3.2.\newline  The Collapsing Conjecture greatly improves our understanding of the Ext--groups in ${\cal P}_d$. For example, it allows  to quickly obtain all computations of \cite{FFSS, C2} (see Cor. 3.8, Cor. 3.9) and reduce those of \cite{C3, C4} to easier particular cases. It has also more conceptual applications
like an intrinsic proof of twist--injectivity in ${\cal P}_d$ (Cor. 3.10). Moreover, applications of the Derived Kan Extension we offer in the present article are not limited to the Collapsing Conjecture and its consequences. We compute in Prop. 4.2 ${\bf K}_{I^{(i)}}^r$ for certain class of Schur functors from which we deduce (Cor. 4.3) all results of \cite{C4} again in a much simpler manner than it was originally done.\newline Beyond obtaining known results in a simpler way, the Derived Kan Extension also helps to understand  general structure of
Ext--groups in ${\cal P}_d$. For example it allows us to  reformulate in Prop. 4.1 the program of computing the Ext--groups between twisted Weyl and Schur functors started in \cite{C2, C4}
in a much more concrete way. This should guide future works on ${\cal P}_d$.
Yet another possible further development of the ideas of our paper would be to investigate ${\bf K}^r_H$  for some other $H$  and apply it to computing Ext--groups for precompositions with
$H$. We briefly discuss this issue at the end of Section 2.\newline
The paper is organized as follows. We define the functors ${\bf K}_H^r, {\bf K}_H^l$ in Section 2. Section 3 contains the proof of the Collapsing Conjecture
and its applications to computations of Ext--groups between twisted functors. In section 4 we obtain further applications to Ext--groups between twisted Weyl
and Schur functors and we discuss  possible further developments.

\section{Precomposition and its adjoints}
Let ${\cal V}$ (resp. ${\cal V}^{gr}$) be the category of finite--dimensional vector spaces (resp. {\bf Z}--graded, totally finite--dimensional, vector
spaces) over a fixed field {\bf k}. We work in the category
 ${\cal P}_{d}$ of homogeneous  strict polynomial functors of degree $d$ over  ${\bf k}$ (see \cite{FS}, Section 2).
 Let ${\cal DP}_d$ denote the bounded derived category of $\pe_d$ i.e.
 it consists of bounded above cohomological complexes  of objects of ${\cal P}_d$ with bounded cohomology modulo quasi--isomorphisms (we recall that $\pe_d$ has a finite (co)homological dimension (c.f. \cite{To}).\newline
For $H\in{\cal P}_s$ let ${\bf C}_H: {\cal P}_d\longrightarrow{\cal P}_{ds}$ be the  functor of taking precomposition with $H$ i.e. ${\bf C}_H(F)(V):=F(H(V))$.
Since ${\bf C}_H$ is an exact functor between abelian categories, it extends degree--wise to their derived categories and we shall, slightly abusing
notation, denote both functors by ${\bf C}_H$.\newline
Our task in this section is to construct functors
\[{\bf K}_H^r, {\bf K}_H^l: {\cal DP}_{ds}\longrightarrow {\cal DP}_{d}\]
which are respectively the right and left adjoint to  ${\bf C}_H$.
We call these functors respectively the right and left (Derived) Kan Extensions by analogy with topology. In fact, the existence of adjoints to the
precomposition is quite general categorical phenomenon. In our context it follows from a suitable version of the Yoneda lemma.\newline
We start with recalling some standard constructions in ${\cal P}_d$.  Let $I^d, D^d, S^d, \Lambda^d\in{\cal P}_d$ denote respectively the tensor, divided, symmetric and exterior power  functors. More generally, let $\mu=(\mu_1,\ldots,\mu_k)$ be a sequence of non--negative integers with sum  $d$. We call such a sequence a partition of $d$ into at most $k$ parts and denote the set of all such sequences by $Q(d,k)$. Then we
put $D^{\mu}:=D^{\mu_1}\otimes\ldots\otimes
D^{\mu_k}$ and we define analogously $S^{\mu},\Lambda^{\mu}$.\newline
We shall also occasionally use strict polynomial bifunctors in the sense of \cite{FF}. We shall denote the category of strict polynomial bifunctors
contravariant of degree $d$ in the first variable and covariant of degree $e$ in the second variable by $\pe_e^d$.
In particular, $\pe^d:=\pe_0^d$ is just the category of contravariant strict polynomial functors of degree $d$ which may be identified with
$(\pe_d)^{op}$.
Similarly, we denote
the category of strict polynomial bifunctors covariant in both variables of degrees $d, e$ by $\pe_{d,e}$.\newline
For any $V\in {\cal V}$ and $F\in\pe_e^d$ we define $F_V\in \pe^d$ and $F^V\in \pe_e$ by the formulae: \[F_V(W):=F(W,V),\] \[F^V(W):=F(V,W).\]
Similarly, for $G\in\pe_{d,e}$ we define $G_V\in\pe_d$ by \[G_V(W):=G(V,W).\]
For $F\in\pe_d$ we define $F_V,F^V\in\pe_d$ by the formulae: \[F_V(W):=F(V\otimes W),\] \[F^V(W):=F(V^*\otimes W).\] Of course, the constructions
for functors and bifunctors are closely related. Namely, for $F\in\pe_d$ we have $F^h\in\pe^d_d, F^t\in\pe_{d,d}$ given by:
\[F^h(V,W):=F(V^*\otimes W),\] \[F^t(V,W):=F(V\otimes W).\] Then we have $F_V=(F^t)_V, F^V=(F^h)^V$.
We shall often refer to functors obtained by these constructions as to ``functors parameterized by $V$''.\newline
We have a left balanced functor \[{\cal H}\!om_{{\cal P}_{e}}: ({\cal P}_{e}^d)^{op}\times {\cal P}_e\ra {\cal P}_d\]
given by the formula \[{\cal H}\!om_{{\cal P}_e}(F,G)(V):=\mbox{Hom}_{{\cal P}_e}(F^V,G).\]
Similarly (and we will use the same notation) we have
\[{\cal H}\!om_{{\cal P}_{e}}: ({\cal P}_{e})^{op}\times {\cal P}_{d,e}\ra {\cal P}_d\]
given by the formula \[{\cal H}\!om_{{\cal P}_e}(F,G)(V):=\mbox{Hom}_{{\cal P}_e}(F,G_V).\]
Let $D^d(I^*\otimes I)\in\pe_{d}^d$ denote the bifunctor
\[(V,W)\mapsto D^d(V^*\otimes W).\]
Then the Yoneda lemma in $\pe_d$ (\cite{FS}, Th. 2.10) can be rephrased as the isomorphism   of functors
\[ {\cal H}\!om_{\pe_d}(D^d(I^*\otimes I),-)\simeq Id.\]
Dually, we define
\[D^d(I\otimes I)(V,W):= D^d(V\otimes W),\]
and we have in the bi-covariant setting  the isomorphism of contravariant functors
\[ {\cal H}\!om_{\pe_d}(-,S^d(I\otimes I))\simeq Id^{\#},\]
where $Id^{\#}:\pe_d^{op}\ra\pe_d$ is the Kuhn duality: $F^{\#}(V):=(F(V^*))^*$.
The more usual way of formulating the Yoneda lemma is by using parameterized functors. Namely, we can  say that we have ``natural in $F$ and $V$
isomorphisms'':
\[\mbox{Hom}_{{\cal P}_{d}}((D^{d})^V,F)\simeq F(V),\]
\[\mbox{Hom}_{{\cal P}_{d}}(F,S^{d}_V)\simeq F^{\#}(V).\]
The Yoneda lemma can be naturally extended to the case when $F$ itself is a bifunctor. One should remember that in that case $F^{\#}$ means dualizing only
with respect to one variable. For example if $F\in {\cal P}_d$ then
\[\mbox{Hom}_{{\cal P}_d}(F^U, S^d_V)\simeq (F^{\#})_U(V).\]
From the Yoneda lemma  one can conclude (\cite{FS}, Th. 2.10) that
 $\{D^d_V\}$ form a family of projective generators of ${\cal P}_d$.
Moreover, since \[D^d_V\simeq\bigoplus_{\mu\in Q(d,k)} D^{\mu_1}\otimes\ldots\otimes D^{\mu_k},\] for $k=\dim(V)$, the family $\{D^{\mu}\}$
for $\mu$ ranging over the set of partitions of $d$ is another set
of projective generators of ${\cal P}_d$.
Analogously,  $\{S^d_V\}$ and $\{S^{\mu}\}$ are  families of injective cogenerators of ${\cal P}_d$.\newline Now we are ready to define our adjoints.
Since ${\cal H}\!om_{\pe_e}$ is  left balanced  we can consider its full derived functor
\[{\cal RH}\!om_{{\cal P}_{e}}: ({\cal DP}_{e}^d)^{op}\times {\cal DP}_e\ra {\cal DP}_d\]
and its bi-covariant counterpart
\[{\cal RH}\!om_{{\cal P}_{e}}: ({\cal DP}_{e})^{op}\times {\cal DP}_{d,e}\ra {\cal DP}_d.\]
In order to ensure that our constructions preserve boundedness of cohomology objects we recall that ${\cal P}_e$ has finite cohomological dimension (in fact
the categories of bifunctors also have finite cohomological dimension but we do not need this stronger fact).\newline
Let, for $H\in\pe_s$, $D^d(I^*\otimes H)\in\pe_{ds}^d$ denote the bifunctor given by the formula
\[(V,W)\mapsto D^d(V^*\otimes H(W)).\]
Then we put
\[{\bf K}_H^r(F)={\cal RH}\!om_{{\cal P}_{ds}}(D^d(I^*\otimes H),F).\]
Using parameterized functors we can rewrite the above formula as
\[{\bf K}^r_H(F)(V)=\mbox{RHom}_{\pe_{ds}}({\bf C}_H((D^d)^V),F).\]

We define the left adjoint dually
\[{\bf K}_H^l(F):=({\bf K}_{H^{\#}}^r(F^{\#}))^{\#}.\]
Since $D^d(I^*\otimes H)\in{\cal P}_{sd}^d$, we get the functors
${\bf K}_H^r, {\bf K}_H^l: {\cal DP}_{ds}\longrightarrow {\cal DP}_{d}$.
\newline
We start with computing the values of ${\bf K}_H^r, {\bf K}_H^l$ on (co)generators.
\begin{prop}
${\bf K}_H^r(S^{ds}_V)=S^d_{H^{\#}(V)}$ and dually
${\bf K}_H^l(D^{ds}_V)=D^d_{H^{\#}(V)}.$
\end{prop}
{\bf Proof:\ }
Since $(D^d)^{\#}\simeq S^d$, we get
\[{\bf K}_H^r(S^{ds}_V)(W)={\cal RH}\!om_{\pe_{ds}}(D^d(I^*\otimes H), S^{ds}_V)(W)\simeq
{\cal H}\!om_{\pe_{ds}}(D^d(I^*\otimes H), S^{ds}_V)(W)\simeq\]
\[{\cal H}\!om_{\pe_{ds}}((D^d(I^*\otimes H))^W, S^{ds}_V)\simeq
{\cal H}\!om_{\pe_{ds}}(({\bf C}_H(D^d))^W, S^{ds}_V)\simeq
(({\bf C}_H(D^d))^{\#})_W(V)\simeq\]
\[({\bf C}_{H^{\#}}((D^d)^{\#}))_W(V)\simeq
({\bf C}_{H^{\#}}(S^d))_W(V)\simeq S^d(H^{\#}(V)\otimes W)\simeq S^d_{H^{\#}(V)}(W).\]

The formula for  ${\bf K}_H^l(D^{ds}_V)$ is proved analogously.\qed
From this computation we deduce our adjunction in a purely formal manner.
\begin{theo}
${\bf K}_H^r,{\bf K}_H^l:{\cal DP}_{ds}\longrightarrow {\cal DP}_{d}$ are respectively
the right and left adjoint functor to ${\bf C}_H:{\cal DP}_{d}\longrightarrow {\cal DP}_{ds}$.
\end{theo}
{\bf Proof:\ }
We only prove the right adjunction, the proof for the left one is analogous. Since $S^{ds}_V$ cogenerate ${\cal P}_{ds}$ it suffices to
establish a natural in $F\in {\cal DP}_d$ and $S^{ds}_V\in \pe_{ds}$ isomorphism
\[\mbox{Hom}_{{\cal DP}_{ds}}({\bf C}_H(F),S^{ds}_V)\simeq \mbox{Hom}_{{\cal DP}_{d}}(F,{\bf K}_H^r(S^{ds}_V)).\]
Applying the Yoneda lemma to the left-hand side we get
\[\mbox{Hom}_{{\cal DP}_{ds}}({\bf C}_H(F), S^{ds}_V)\simeq
\mbox{Hom}_{\pe_{ds}}({\bf C}_H(F), S^{ds}_V)\simeq ({\bf C}_H(F))^{\#}(V)\simeq F^{\#}(H^{\#}(V)).\]
Using additionally Prop. 2.1 we get at the right-hand side
\[\mbox{Hom}_{{\cal DP}_{d}}(F,{\bf K}_H^r(S^{ds}_{V}))\simeq \mbox{Hom}_{{\cal DP}_{d}}(F,S^d_{H^{\#}}(V))\simeq (F^{\#}(H^{\#}(V)),\]
which proves our statement.\qed
{\bf Remark: } Of course, the same constructions produce the two--sided adjunction ${\cal P}_d\rec{\cal P}_{sd}$ already at the level of abelian categories.
In fact this adjunction for $H=I^{(i)}$ (the $i$th Frobenius twist functor) can be found among those considered by Kuhn in \cite{Ku}. We will extend this
remark in Section 3.\newline
Let us now consider  the problem of computing the groups $\mbox{Ext}^*_{{\cal P}_{ds}}({\bf C}_H(F),G)$ for $F\in{\cal P}_d, G\in {\cal P}_{ds}$.
By our adjunction, once we compute ${\bf K}^r_H(G)$, the calculations can be transferred  to ${\cal P}_d$. On the other hand, since
$\bigoplus (D^d)^V\simeq D^{\mu_1}\otimes\ldots\otimes D^{\mu_k}$ (non--naturally in $V$),
computing ${\bf K}^r_H(G)$ requires, roughly speaking,  computing $\mbox{Ext}^*_{{\cal P}_{ds}}({\bf C}_H(D^{\mu}),G)$, organizing results functorially
in $V$ and understanding the hyper--Ext spectral sequence.  Thus a general problem is reduced to the special case of $F=D^{\mu}$ and computations in smaller degrees. We will see in the next sections
how this procedure works for $H=I^{(i)}$ (the i--fold Frobenius twist functor) but it would also be interesting to apply this approach to precompositions
with other $H$ like $D^s, S^s, \Lambda^s$ etc. \newline We finish this section with considering certain general condition on $F\in\pe_d$ under which   $H^*({\bf K}^r_H(F))$
can be described quite explicitly. This criterion will be vital
 in the applications discussed in Section 4.
Let $\Sigma_{\mu}=\Sigma_{{\mu}_1}\times\ldots\Sigma_{{\mu}_k}\subseteq\Sigma_d$
be   the Young subgroup associated to a partition $\mu=\{\mu_1,\ldots,\mu_k\}$. Then  $D^{\mu}\simeq(I^d)^{\Sigma_{\mu}}.$
\begin{prop}
Assume that for any Young subgroup $\Sigma_{\mu}\subseteq\Sigma_d$,
the natural embedding $D^{\mu}\longrightarrow I^d$ induces an isomorphism on Ext--groups
\[\mbox{Ext}^*_{{\cal P}_{ds}}({\bf C}_H(D^{\mu}),F)\simeq (\mbox{Ext}^*_{{\cal P}_{ds}}({\bf C}_H(I^d),F))_{\Sigma_{\mu}}.\]
Then
\[H^*({\bf K}_H^r(F))\simeq I^d\otimes_{\Sigma_d} \mbox{Ext}^*_{{\cal P}_{ds}}({\bf C}_H(I^d),F).\]
Dually, if \[\mbox{Ext}^*_{{\cal P}_{ds}}(F,{\bf C}_H(S^{\mu}))\simeq (\mbox{Ext}^*_{{\cal P}_{ds}}(F,{\bf C}_H(I^d)))_{\Sigma_{\mu}},\]
then
\[H^*({\bf K}_H^l(F))\simeq\mbox{Hom}_{\Sigma_d}(\mbox{Ext}^*_{{\cal P}_{ds}}(F,{\bf C}_H(I^d)), I^{d}).\]
\end{prop} {\bf Proof:\ } Again we  focus on the formula for ${\bf K}_H^r$, leaving the case of ${\bf K}_H^l$ to the reader.
We consider the composition
\[\alpha: \mbox{Hom}_{{\cal P}_{d}}((D^d)^V,I^d)\otimes \mbox{RHom}_{{\cal P}_{ds}}({\bf C}_H(I^d), F)\longrightarrow\]
\[\mbox{Hom}_{{\cal P}_{ds}}({\bf C}_H((D^d)^V),{\bf C}_H(I^d))\otimes \mbox{RHom}_{{\cal P}_{ds}}({\bf C}_H(I^d),F)\longrightarrow \mbox{RHom}_{{\cal P}_{ds}}({\bf C}((D^d)^V),F)\]
where the first arrow is the precomposition with $H$ tensored with the identity and the second is the Yoneda composition.
We shall look at the effect of this transformation on the Ext--groups (i.e. cohomology of RHom).
We take $\mu\in Q(d,n)$ for $n=\dim(V)$ and  choose some basis of $V$.
Since by the Yoneda Lemma, $\mbox{Hom}_{{\cal P}_d}((D^d)^V,I^d)\simeq V^{\otimes d}$,
the subspace of $\mbox{Hom}_{{\cal P}_d}((D^d)^V,I^d)\otimes \mbox{Ext}^*_{{\cal P}_{ds}}({\bf C}_H(I^d), F)$ on which
the standard torus in $GL_n(\ka)$ acts with  weight $\mu$ may be identified with
$({\bf k}\otimes_{\Sigma_{\mu}}{\bf k}[\Sigma_d])\otimes \mbox{Ext}^*_{{\cal P}_{ds}}({\bf C}_H(I^d),F)$.
On the other hand the $\mu$-weight subspace of $\mbox{Ext}^*_{{\cal P}{ds}}({\bf C}_{H}((D^d)^V),F)$ may be identified with
$\mbox{Ext}^*_{{\cal P}_{ds}}({\bf C}_H(D^{\mu}),F)$. Under these identifications  $\alpha_*$ is just
induced by the embedding $D^{\mu}\longrightarrow I^{d}$. Thus, by our assumption  $\alpha_*$  factorizes to the isomorphism
\[({\bf k}\otimes_{\Sigma_{\mu}}{\bf k}[\Sigma_d])\otimes_{\Sigma_d} \mbox{Ext}^*_{{\cal P}_{ds}}({\bf C}_H(I^d),F)\simeq
\mbox{Ext}^*_{{\cal P}_{ds}}({\bf C}_H(D^{\mu}),F).\]
Then by gathering up these isomorphisms for all $\mu$ we get an isomorphism
\[V^{\otimes d}\otimes_{\Sigma_d} \mbox{Ext}^*({\bf C}_H(I^d),F)\simeq \mbox{Ext}^*_{{\cal P}_{ds}}({\bf C}_H((D^d)^V),F)\simeq H^*({\bf K}_H^r(F))(V)
.\]\qed

\section{The Collapsing Conjecture and related Ext--groups}
In the rest of the paper we restrict our  attention to the case of a base field ${\bf k}$ of positive characteristic $p$ and $H=I^{(i)}$ (the i-th Frobenius twist). Hence from now on we denote ${\bf C}_{I^{(i)}}$ just by ${\bf C}_i$,
${\bf K}_{I^{(i)}}^r$  by ${\bf K}_i$ and ${\bf K}^r_{I^{(i)}}\circ {\bf C}_{I^{(i)}}$ by ${\bf KC}_i$ (still all our results  have obvious counterparts for
${\bf K}_{I^{(i)}}^l$ but we leave their formulation to the reader). In order  to further simplify notation we rewrite the formula defining ${\bf K}_i$ as
\[{\bf K}_i(F)(V)=\mbox{RHom}_{\pe_{pd}}(D^{d(1)}_{V^*},F),\]
although one should remember that in fact $D^{d(1)}_{V^*}={\bf C}_i(D^d_{V^*})$ i.e. we do not twist the parameter $V$. Thus we have a somewhat unfortunate
fact that $D^{d(1)}_{V^*}\neq (D^{d(1)})_{V^*}$ but we will not use the latter functor in our article.
We also prefer here and later in several places the notation $F_{V^*}$ instead of $F^V$ to avoid overloading superscript.\newline
In the present section we prove a stronger version of Touz\'e's Collapsing Conjecture
and derive some its consequences concerning the groups $\mbox{Ext}^*_{{\cal P}_{dp^i}}(F^{(i)}, G^{(i)})$.
 We start  with  a simple observation which makes use of multiplicativity of the Frobenius twist. We recall that if we apply a strict polynomial functor $F$
 to a graded space  $V$ then $F(V)$ has a natural grading (see e.g. \cite{T2}, Sect.~2.5).
 This grading is self--evident for functors of tensor type but the case of Frobenius twist functor is a bit tricky: one should remember that
 $(V[j])^{(1)}:=V^{(1)}[pj]$.\newline
 Now, for $F\in\pe_d$ and $V\in{\cal V}^{gr}$, $F_V$ is
 a graded strict polynomial functor hence it may be regarded as a complex with trivial differential.
 Moreover, if we have $F,G\in {\cal P}_d$ and $\alpha: F\ra G$, then $\alpha_{V}:F_{V}\ra G_{V}$ preserves grading. Therefore, if $C^{\bullet}$ is a complex of objects
of ${\cal P}_d$ then  $C^{\bullet}_{V}$ may be naturally thought of as a complex, hence an object of ${\cal DP}_d$.
\begin{prop}
For any $V\in{\cal V}^{gr}$, and $F\in {\cal DP}_{dp^i}$,      ${\bf K}_i(F_V)\simeq ({\bf K}_i(F))_{V^{(i)}}$.
\end{prop}
{\bf Proof:\ } It suffices to establish this formula for  $F=S^{ds}_W$.
We recall that $(I^{(i)})^{\#}\simeq I^{(i)}$. Then by Proposition 2.1
\[{\bf K}_i((S^{dp^i}_W)_V)={\bf K}_i(S^{dp^i}_{W\otimes V})\simeq S^d_{(W\otimes V)^{(i)}}\simeq S^d_{W^{(i)}\otimes V^{(i)}}\simeq (S^d_{W^{(i)}})_{V^{(i)}}\simeq \]
\[\simeq({\bf K}_i(S^{dp^i}_{W}))_{V^{(i)}}.\]
\qed
We now turn to the main objective of this section i.e. computing ${\bf KC}_i(F)={\bf K}_i(F^{(i)})$ for any $F\in{\cal DP}_d$.
Let $A:={\bf k}[x]/x^p$ for $|x|=2$, and more generally $A_i:=A\ot A^{(1)}\ot\ldots\ot A^{(i-1)}$. Then by \cite{FS}, Th.~4.5,
${\mbox Ext}^*_{{\cal P}_{p^i}}(I^{(i)},I^{(i)})$ as a graded {\bf k}-algebra with Yoneda multiplication is isomorphic to $A_i$. \newline
We recall that the assignment $F\mapsto F_{A_i}$  produces an exact endofunctor on the category ${\cal K}\!om({\cal P}_d)$.
 Hence it gives the functor $Id_{A_i}: {\cal DP}_d\ra {\cal DP}_d$.
\begin{theo}
There is an isomorphism of functors $Id_{A_i}\simeq {\bf KC}_i$.
\end{theo}
We need for the proof of Theorem 3.2  certain elements  $e_d\in \mbox{Hom}_{{\cal DP}^d_{pd}}(D^d(I^*\otimes I^{(1)}),D^d(I^*\otimes I^{(1)})_{A^*})$
where $D^d(I^*\otimes I^{(1)})_{A^*}(V,W):=D^d(V^*\ot V^{(1)}\ot A^*)$ and $A^*$ is the graded (hence concentrated in non--positive degrees) linear dual
of $A$. We start by constructing these classes. Let us first look at the case $d=1$. Then we have (by using the K\"unneth formula (\cite{FFSS}, p.~672) for the third
isomorphism)
\[\mbox{Hom}_{{\cal DP}^1_{p}}(I^*\otimes I^{(1)},(I^*\otimes I^{(1)})_{A^*})\simeq
\mbox{Hom}_{{\cal DP}^1_{p}}(I^*\otimes I^{(1)},I^*\otimes I^{(1)}\otimes A^*)\simeq\] \[\bigoplus_{j=0}^{p-1}
\mbox{Hom}_{{\cal DP}^1_{p}}(I^*\otimes I^{(1)},I^*\otimes I^{(1)}[2j])\simeq
\bigoplus_{j=0}^{p-1}\mbox{Hom}_{{\cal P}^1}(I^*,I^*)\otimes \mbox{Ext}^{2j}_{{\cal P}_p}(I^{(1)},I^{(1)})\simeq\]\[
\bigoplus_{j=0}^{p-1}\mbox{Ext}^{2j}_{{\cal P}_p}(I^{(1)},I^{(1)})\simeq A.\]
Under this identification we put $e_1$ to be $\oplus_{j=0}^{p-1}x^j$ and we claim that we can also choose similar elements
for higher $d$ in a compatible manner.
\begin{prop}
There exist classes \[e_d\in \mbox{Hom}_{{\cal DP}^d_{pd}}(D^d(I^*\otimes I^{(1)}),(D^d(I^*\otimes I^{(1)}))_{A^*})\] satisfying:
\begin{enumerate}
\item $e_1=\oplus_{j=0}^{p-1}x^j$.
\item $e_1^{\ot d}\circ \Delta=\Delta_{A*}\circ e_d$ as elements of
\[\mbox{Hom}_{{\cal DP}^d_{pd}}(D^d(I^*\otimes I^{(1)}),I^d(I^*\otimes I^{(1)})_{A^*}),\]
where $\Delta: D^d(I^*\ot I)\ra I^{d}(I^*\ot I)$ is the natural embedding and
\[e_1^{\ot d}\in \mbox{Hom}_{{\cal DP}^d_{pd}}(I^d(I^*\otimes I^{(1)}),I^d(I^*\otimes I^{(1)})_{A^*})\]
is the $d$th external power of  $e_1$ (see e.g. \cite{T1}, Sect.~1.1).
\end{enumerate}
\end{prop}
We will obtain our classes $e_d$ by gathering up some elements in Ext--groups which are closely
related to the Touz\'e classes \cite{T1,T3}.
The following sub--lemma will allow us to use the Touz\'e classes.
\begin{lem}
For any partition $\mu$ of $d$
there is an isomorphism
\[\mbox{Ext}^*_{{\cal P}^d_{pd}}(D^{d}(I^*\otimes I^{(1)}),D^{\mu}(I^*\otimes I^{(1)}))\simeq
\mbox{Ext}^*_{{\cal P}^{pd}_{pd}}(D^{pd}(I^*\otimes I),D^{\mu(1)}(I^*\otimes I))\]
natural with respect to maps between  $D^{\mu}$.
\end{lem}
This lemma follows from the existence of spectral sequence which is an incarnation of the classical ``associativity formula'' (\cite{CE}, XVI. 4. (5)) and may be of some independent interest.
\begin{prop}
For any $F,G\in{\cal P}^d_e$ there exists a spectral sequence
\[E_2^{st}=\mbox{Ext}^s_{{\cal P}_e^e}(D^e(I^*\ot I),\mbox{Ext}_{{\cal P}^d}^t(F,G)) \Longrightarrow \mbox{Ext}^{s+t}_{{\cal P}_e^d}(F,G)\]
This spectral spectral sequence is natural in $F,G$ and commutes with the external product. By the last assertion we mean the following.\\
For any $F_i,G_i\in {\cal P}_{e_i}^{d_i},\ i=1,2,$ we have a map of spectral sequences
\[E(F_1,G_1)\ot E(F_2,G_2)\ra E(F_1\ot F_2,G_1\ot G_2).\]
We require that:
\begin{itemize}
\item The map  at the $E_2$ pages coincides with
\[\mbox{Ext}^{s_1}_{{\cal P}_{e_1}^{e_1}}(D^{e_1}(I^*\ot I),\mbox{Ext}_{{\cal P}^{d_1}}^{t_1}(F_1,G_1))\ot
\mbox{Ext}^{s_2}_{{\cal P}_{e_2}^{e_2}}(D^{e_2}(I^*\ot I),\mbox{Ext}_{{\cal P}^{d_2}}^{t_2}(F_2,G_2))\ra\]
\[\mbox{Ext}^{s_1+s_2}_{{\cal P}_{e_1+e_2}^{e_1+e_2}}(D^{e_1}(I^*\ot I)\ot D^{e_2}(I^*\ot I),
\mbox{Ext}_{{\cal P}^{d_1}}^{t_1}(F_1,G_1)\ot \mbox{Ext}_{{\cal P}^{d_2}}^{t_1}(F_2,G_2))\ra\]
\[\mbox{Ext}^{s_1+s_2}_{{\cal P}_{e_1+e_2}^{e_1+e_2}}(D^{e_1}(I^*\ot I)\ot D^{e_2}(I^*\ot I),
\mbox{Ext}_{{\cal P}^{d_1+d_2}}^{t_1+t_2}(F_1\ot F_2,G_1\ot G_2))\ra\]
\[\mbox{Ext}^{s_1+s_2}_{{\cal P}_{e_1+e_2}^{e_1+e_2}}(D^{e_1+e_2}(I^*\ot I),
\mbox{Ext}_{{\cal P}^{d_1+d_2}}^{t_1+t_2}(F_1\ot F_2,G_1\ot G_2))\]
where the first two arrows are external products and the third is induced by the embedding $D^{e_1+e_2}\ra D^{e_1}\ot D^{e_2}$.
\item The map on the abutments coincides (up to filtration) with the external product
\[\mbox{Ext}_{{\cal P}_{e_1}^{d_1}}^{s_1}(F_1,G_1)\ot \mbox{Ext}_{{\cal P}_{e_2}^{d_2}}^{s_2}(F_2,G_2))\ra
\mbox{Ext}_{{\cal P}^{d_1+d_2}_{e_1+e_2}}^{s_1+s_2}(F_1\ot F_2,G_1\ot G_2)).\]
\end{itemize}
\end{prop}
{\bf Proof of Proposition 3.5: } We consider the functors
\[\alpha_F: {\cal P}_e^d\ra {\cal V},\ \  \alpha_F(G):=\mbox{Hom}_{{\cal P}_d^e}(F,G),\]
\[\beta_F:{\cal P}_e^d\ra{\cal P}_e^e,\ \  \beta_F(G):={\cal H}\!om_{{\cal P}^d}(F,G),\]
\[\gamma:{{\cal P}_e^e}\ra{\cal V},\ \ \gamma(G):=\mbox{Hom}_{{\cal P}_e^e}(D^e(I^*\ot I),G).\]
We have then $\alpha_F\simeq\gamma\circ\beta_F$. To see this, it suffices to evaluate both sides on injective cogenerators
of ${\cal P}_e^d$ which are given by $I_{X,Y}^{d,e}:=(D^d_{X^*})^*\ot S^e_{Y^*}$  (\cite{FF}, Prop.~1.2).
Then the suitable version of the Yoneda lemma (see the proof of \cite{FF}, Prop.~1.2) says that
\[\mbox{Hom}_{{\cal P}_e^d}(F,I^{d,e}_{X,Y})\simeq F^*(X,Y)\]
for any $F\in{\cal P}_e^d$ ($F^*$ for a (bi)functor $F$ stands for $F$ postcomposed with $(-)^*$).
Therefore we obtain
\[\beta_F(I^{e,d}_{X,Y})={\cal H}\!om_{{\cal P}^d}(F,(D^d_{X^*})^*\ot S^e_{Y^*})\simeq
{\cal H}\!om_{{\cal P}^d}(F,(D^d_{X^*})^*)\ot S^e_{Y^*}\simeq\]
\[{\cal H}\!om_{{\cal P}_d}(D^d_{X^*},F^*)\ot S^e_{Y^*}\simeq F^*(X,-)\ot S^e_{Y^*}.\]
Then by using (\cite{FF}, Prop.~1.3) which says that
\[\mbox{Hom}_{{\cal P}_e^e}(D^e(I^*\ot I),F^*\ot G)\simeq \mbox{Hom}_{{\cal P}_e}(F,G)\]
for any $F,G\in{\cal P}_e$, we conclude that
\[\gamma\circ\beta_F(I^{e,d}_{X,Y})=\mbox{Hom}_{{\cal P}_e^e}(D^e(I^*\ot I),F^*(X,-)\ot S^e_{Y^*})\simeq
\mbox{Hom}_{{\cal P}_e}(F(X,-),S^e_{Y^*})\simeq\]\[ F^*(X,Y)\simeq\alpha_F(I^{e,d}_{X,Y}).\]
Moreover, by \cite{FF}, Prop.~1.3 again, for any $F\in{\cal P}_e$ the bifunctor $F^*(X,-)\ot S^e_{Y^*}$ is  $\gamma$-acyclic, hence $\beta_F$ takes injective generators to
$\gamma$--acyclic objects. Thus we get our spectral sequence
as the Grothendieck spectral sequence associated to the composite $\gamma\circ \beta_F$. The naturality and commuting with external product of our spectral
sequence follows from the naturality of the Grothendieck spectral sequence and of the external product.\qed
{\bf Proof of Lemma 3.4: } We consider the  spectral sequence from Proposition 3.5 for $F=D^{d}(I^*\ot I^{(1)}), G=D^{\mu}(I^*\ot I^{(1)})$.
Then, since $F$ is projective in the contravariant variable, the spectral sequence collapses at $E_2$ and we get our assertion.\qed
{\bf Proof of Proposition 3.3: }
Let $c_d\in\mbox{Ext}^{2d}_{{\cal P}^d_{pd}}(D^{d}(I^*\otimes I^{(1)}),D^{d}(I^*\otimes I^{(1)}))$ be the dth Touz\'e class (see \cite{T1}, Th.~1.4) pulled by the isomorphism of Lemma 3.4 for $\mu=(d)$ and let $c_d^j\in\mbox{Ext}^{2dj}_{{\cal P}^d_{pd}}(D^{d}(I^*\otimes I^{(1)}),D^{d}(I^*\otimes I^{(1)}))$ be its $j$th Yoneda power.
Now we have
\[(D^d(I^*\ot I^{(1)}))_{A^*}\simeq\bigoplus_{\mu\in Q(d,p)}\bigotimes_j D^{\mu_j}(I^*\ot I^{(1)})[2(j-1)\mu_j].\]
Let $\Delta_{\mu}$ be the embedding $D^{d}(I^*\ot I^{(1)})\ra D^{\mu}(I^*\ot I^{(1)})$. Then we put $e_d:=\bigoplus_{\mu}e_{\mu}$ where $e_{\mu}
\in \mbox{Ext}^{\Sigma\, 2(j-1)\mu_j}_{{\cal P}_{pd}^d}(D^d(I^*\ot I^{(1)}),D^{\mu}(I^*\ot I^{(1)}))$ equals
$(\bigotimes_j c_{\mu_j}^{j-1})\circ \Delta_{\mu}$. Then it is easy to see that $e_1=\oplus_{j=0}^{p-1}x^j$. \newline
To establish the second condition of Proposition 3.3
we start by considering the pulled Touz\'e classes $c_d\in\mbox{Ext}^{2d}_{{\cal P}^d_{pd}}(D^{d}(I^*\otimes I^{(1)}),D^{d}(I^*\otimes I^{(1)}))$
again. Let  $\widetilde{c_d}\in\mbox{Ext}^{2d}_{{\cal P}^{pd}_{pd}}(D^{pd}(I^*\otimes I),D^{d(1)}(I^*\otimes I))$ be the original Touz\'e class
and let
\[\chi_{\mu}:
\mbox{Ext}^*_{{\cal P}^{pd}_{pd}}(D^{pd}(I^*\otimes I),D^{\mu(1)}(I^*\otimes I))\ra
\mbox{Ext}^*_{{\cal P}^d_{pd}}(D^{d}(I^*\otimes I^{(1)}),D^{\mu}(I^*\otimes I^{(1)}))\]
be the isomorphism from Lemma 3.4. Then thanks to the naturality and commuting with external product of the spectral sequence producing $\chi_{\mu}$
we have
\[c_1^{\ot d}\circ\Delta=(\chi_{(1)}(\widetilde{c_1}))^{\ot d}\circ\Delta=\chi_{(1^d)}(\widetilde{c_1}^{\ot d}\circ\Delta).\]
Next, using the naturality of $\chi_{\mu}$ again and [T1,Th.~1.4] we obtain
\[\Delta\circ c_d=\Delta\circ \chi_{(d)}(\widetilde{c_d})=\chi_{(1^d)}(\Delta\circ \widetilde{c_d})=\chi_{(1^d)}(\widetilde{c_1}^{\ot d}\circ \Delta)=
c_1^{\ot d}\circ\Delta.\]

Then we get immediately the analogous relation \[\Delta\circ c_d^j=(c_1^j)^{\ot d}\circ \Delta\] for the Yoneda powers of $c_d$.\newline
We recall from the definition of $e_d$ that   $D^d(I^*\ot I^{(1)}))_{A^*}$ can be decomposed into $\bigoplus_{\mu} D^{\mu}(I^*\ot I^{(1)})[|2(j-1)\mu|]$ (where $|2(j-1)\mu|$ stands for $\sum_j 2(j-1)\mu_j$) and, with this decomposition, $e_d=\bigoplus e_{\mu}$ for $e_{\mu}=c_{\mu}\circ\Delta_{\mu}$ (where
$c_{\mu}:=\bigotimes_j c_{\mu_j}^{j-1}$). We have an analogous decomposition for the tensor power
\[(I^d(I^*\ot I^{(1)}))_{A^*}\simeq\bigoplus_{\mu\in Q(d,p)} (I^{d}(I^*\ot I^{(1)}))\!\!\uparrow_{\Sigma_{\mu}}^{\Sigma_d}[|2(j-1)\mu|].\]
where $(I^{d}(I^*\ot I^{(1)}))\!\!\uparrow_{\Sigma_{\mu}}^{\Sigma_d}=\mbox{Hom}_{{\bf k}[\Sigma_{\mu}]}(\ka[\Sigma_{d}],I^{d}(I^*\ot I^{(1)}))$
is the coinduced $\Sigma_d$--module.
 With this decomposition, $e_1^{\ot d}=\bigoplus_{\mu} e_1^{\mu}$ and $e_1^{\mu}=c_1^{\mu}\!\!\uparrow_{\Sigma_{\mu}}^{\Sigma_d}\circ \Delta_{\lambda/\mu}:=(\bigotimes_j (c_1^{j-1})^{\ot\mu_j})\!\!\uparrow_{\Sigma_{\mu}}^{\Sigma_d}\circ \Delta_{\lambda/\mu}$
where $\Delta_{\lambda/\mu}: I^d(I^*\ot I^{(1)})\ra(I^* \ot I^{(1)})\!\!\uparrow_{\Sigma_{\mu}}^{\Sigma_d}$ is the natural $\Sigma_d$--equivariant
embedding.\newline Let us look at the diagram

\[
\begin{array}{ccccc}
D^d(I^*\ot I^{(1)})&\stackrel{\Delta_{\mu}}{\ra}&D^{\mu}(I^*\ot I^{(1)})&\stackrel{c_{\mu}}{\ra}&D^{\mu}(I^*\ot I^{(1)})[|2(j-1)\mu|]\\
\downarrow \scriptstyle{\Delta}&\mbox{}&\downarrow \scriptstyle{\Delta _{\hat{\mu}}}&\mbox{}&\downarrow \scriptstyle{\Delta_{\hat{\mu}}[|2(j-1)\mu|]}\\
I^d(I^*\ot I^{(1)})&\Delta_{\lambda/\mu}&
(I^{d}(I^*\ot I^{(1)}))\!\!\uparrow_{\Sigma_{\mu}}^{\Sigma_d}&\stackrel{c_1^{\mu}\!\uparrow_{\Sigma_{\mu}}^{\Sigma_d}}{\ra}&
(I^{d}(I^*\ot I^{(1)}))\!\!\uparrow_{\Sigma_{\mu}}^{\Sigma_d}[|2(j-1)\mu|]
\end{array}
\]
where $\Delta_{\hat{\mu}}: D^{\mu}(I^*\ot I^{(1)})\ra (I^d(I^*\ot I^{(1)}))\!\!\uparrow_{\Sigma_{\mu}}^{\Sigma_d}$ is the natural embedding of the $\Sigma_d$--invariants. We claim that this diagram is commutative. Indeed, the commutativity of the left square is obvious, while the commutativity of
of the right one follows from the relation $\Delta\circ c_d^j=(c_1^j)^{\ot d}\circ \Delta$. Now the commutativity of the whole diagram gives \[\Delta\circ e_{\mu}=e_1^{\mu}\circ\Delta.\]
Since $e_d=\bigoplus_{\mu} e_{\mu}$ and $e_1^{\ot d}=\bigoplus_{\mu} e_1^{\mu}$ we get
$ \Delta\circ e_d=e_1^{\ot d}\circ \Delta$ which concludes the proof.
\qed
{\bf Proof of Theorem 3.2: }
We start with computing dimensions of cohomology groups of ${\bf KC}_i$ for injective cogenerators of ${\cal P}_d$.
\begin{lem}
Let $\mu$ be a partition of $d$. Then there is a degree--wise equality
\[\dim (S^{\mu}_{A_i})=\dim (H^*({\bf KC}_i(S^{\mu})).\]
\end{lem}
{\bf Proof of Lemma 3.6: } We first consider $\mu=(1,\ldots,1)$. Then by the K\"unneth formula (\cite{FFSS}, p. 672) we get
\[H^*({\bf KC}_i(I^d))(V)=\mbox{Ext}^*_{\pe_{dp^i}}(D^{d(i)}_{V^*},I^{d(i)})\simeq
\bigotimes_{s=1}^d\mbox{Ext}^*_{\pe_{p^i}}(I^{(i)}_{V^*}, I^{(i)})\simeq\]
\[\simeq\bigotimes_{s=1}^d A_i\otimes V\simeq (I^d)_{A_i}.\]
Now we take an arbitrary $\mu$. Since
\[(D^d)_{V^*}\simeq\bigoplus_{\mu} D^{\rho_1}\otimes\ldots\otimes D^{\rho_k}\] (non-canonically in $V$) and by \cite{C2}, Prop. 4.1
\[\mbox{Ext}^*_{\pe_{dp^i}}(D^{\rho(i)},S^{\mu(i)})\simeq (\mbox{Ext}^*_{\pe_{dp^i}}(D^{\rho(i)},I^{d(i)}))_{\Sigma_{\mu}},\] we get
\[\dim(H^*({\bf KC}_i(S^{\mu})))\simeq \dim (H^*({\bf KC}_i(I^d)) _{\Sigma_{\mu}})\simeq
\dim((I^d_{A_i})_{\Sigma_{\mu}})\simeq \dim(S^{\mu}_{A_i}).\]
\qed
We have now collected all ingredients we need for the construction of isomorphisms $\Phi_d^i: Id_{A_i}\ra {\bf KC}_i$ between endofunctors
on ${\cal DP}_d$. We proceed by induction on $i$. We start with $i=1$.\newline
 Let \[\phi_d:
Id_A\ra {\cal RH}\!om_{{\cal P}_{pd}}(D^d(I^*\ot A^*\ot I^{(1)}), F^{(1)})\simeq {\bf KC}_1\circ Id_A\] be the composite
\[F_A\simeq {\cal H}\!om_{{\cal P}_d}(D^d(I^*\ot A^*\ot I), F)\ra  {\cal H}\!om_{{\cal P}_{pd}}(D^d(I^*\ot A^*\ot I^{(1)}), F^{(1)})
\ra\]\[ {\cal RH}\!om_{{\cal P}_{pd}}(D^d(I^*\ot A^*\ot I^{(1)}), F^{(1)})\]
of the Yoneda isomorphism with twisting and embedding of ${\cal H}\!om$ into ${\cal RH}\!om$  (in fact $\phi_d$ is nothing but the unit of our adjunction
applied to $F_A$ but we will not use this).
Let  $e_d^*: {\bf KC}_1\circ Id_A\ra {\bf KC}_1$ be the transformation induced by $e_d$.
Then we put
\[\Phi_d^1(F):=e_d^*\circ\phi_d(F).\]
We start by showing that $\Phi_1^1(I)$ is a quasi--isomorphism. To this end we shall describe $\Phi_1^1(I)(V)(v\ot a)$ for $v\ot a\in V\ot A$
explicitly as an element of $\mbox{RHom}_{{\cal P}_{p}}(I^{(1)}_{V^*}, I^{(1)})$.
For $v\ot a\in V\ot A$ let $\rho_{v\ot a}: I_{(V\ot A)^*}\ra I$ be the evaluation map $\rho_{v\ot a}(u\ot \alpha):=\alpha(v\ot a)u$. Then the Yoneda isomorphism
$I_A\simeq \mbox{Hom}_{{\cal P}_1}(I_{V\ot A}, I)$ just sends $v\ot a$ to $\rho_{v\ot a}$. Thus $\phi(V)(v\ot a)$ is $\rho_{v\ot a}$ precomposed with $I^{(1)}$
which we will slightly abusing notation also denote by $\rho_{v\ot a}$. Then finally
$\Phi_1^1(I)(V)(v\ot a)=\rho_{v\ot a}\circ e_1$. Now when we identify $\mbox{RHom}_{{\cal P}_{p}}(I^{(1)}_{V^*}, I^{(1)})$ with
$V\ot \mbox{RHom}_{{\cal P}_{p}}(I^{(1)}, I^{(1)})$ we get by Proposition 3.3.1 that $H^*(\rho_{v\ot x^j}\circ e_1)$ corresponds to $v\ot x^j$. This shows that
$\Phi_1^1(I)$ is a quasi--isomorphism as we claimed.\newline
Now we are going to show that  $\Phi_d^1(I^d)$ is a quasi-isomorphism. Again we start by identifying $\Phi_d^1(I^{\ot d})(V)({\bf v}\ot {\bf a}) $
for ${\bf v}\ot {\bf a}\in V^{\ot d}\ot A^{\ot d}$. This time the Yoneda isomorphism sends ${\bf v}\ot {\bf a}$ to
$\rho^d_{{\bf v}\ot {\bf a}}\circ \Delta$ where $\Delta: D^d_{(V\ot A)^*}\ra I^d_{(V\ot A)^*}$ is the embedding and $\rho^d_{{\bf v}\ot {\bf a}}: I^d_{(V\ot A)^*}\ra I^d$ is the evaluation map. Hence we get $\Phi_d^1(I^d)_V({\bf v}\ot {\bf a})=\rho^d_{{\bf v}\ot {\bf a}}\circ \Delta\circ e_d$.
Now by the K\"unneth formula  we have $\mbox{RHom}_{{\cal P}_{pd}}(D^{d(1)}_{V^*}, I^{d(1)})\simeq \bigotimes_{j=1}^d\mbox{RHom}_{{\cal P}_{p}}(I^{(1)}_{V^*}, I^{(1)})$. With this identification, since by Proposition 3.3.2 $\Delta\circ e_d=e_1^{\ot d}\circ\Delta$, we get $H^*(\Phi_d^1(I^d))(V)= \bigotimes_{j=1}^d H^*(\Phi_1^1(I))(V)$. Thus
$\Phi_d^1(I^d)$ is a quasi--isomorphism.\newline
Now as we remember from Lemma 3.6 and its proof, the multiplication map $m_{\lambda}:I^d\ra S^{\lambda}$ induces an epimorphism $m_*: H^*({\bf KC}_1(I^d))\ra H^*({\bf KC}_1(S^{\lambda}))$ and $\mbox{dim}(S^{\lambda}_A)=
\mbox{dim}(H^*({\bf KC}_1(S^{\lambda})))$. Hence we get that $\Phi_d^1(S^{\lambda})$ is a quasi--isomorphism. Since $S^{\lambda}$ cogenerate ${\cal P}_d$ and $\Phi_d^1$
is additive, we conclude that $\Phi_d^1(F)$ is a quasi-isomorphism for any $F\in{\cal DP}_d$.\newline
This finishes the proof of Theorem 3.2 for $i=1$. We now turn to the general case. We assume that we have  isomorphisms
$\Phi_d^i: Id_{A_i}\simeq {\bf KC}_i$ for all $d$. Then since ${\bf C}_{i+1}\simeq {\bf C}_i\circ {\bf C}_1$ and
${\bf K}_{i+1}\simeq {\bf K}_1\circ {\bf K}_i$, we get \[{\bf KC}_{i+1}\simeq
{\bf K}_1\circ {\bf K}_i\circ {\bf C}_i\circ {\bf C}_1\simeq
{\bf K}_1\circ Id_{A_i}\circ {\bf C}_1.\]
Now by Proposition 3.1, $Id_{A_i}\circ{\bf C}_1\simeq {\bf C}_1\circ Id_{A_i^{(1)}}$. Thus
\[{\bf K}_1\circ Id_{A_i}\circ {\bf C}_1\simeq {\bf K}_1\circ {\bf C}_1\circ Id_{A_i^{(1)}}
\simeq (Id_{A_i^{(1)}})_{A}\simeq Id_{A_{i+1}}.\]
\qed

Theorem 3.2 by taking cohomology, gives
\begin{cor}
For any $F,G\in{\cal P}_d$, there is a natural in $F,G$ isomorphism  of graded spaces
\[\mbox{Ext}^*_{{\cal P}_{dp^i}}(F^{(i)},G^{(i)})\simeq \mbox{Ext}^*_{{\cal P}_d}(F,G_{A_i}).\]
\end{cor}

{\bf Remark: } As we have mentioned in Section 2, the abelian part of our adjunction was investigated by Kuhn \cite{Ku}. In particular Kuhn (\cite{Ku}, Th.~6.10) shows
that the two--sided adjunction ${\cal P}_d\rec {\cal P}_{p^i d}$ is a part of recollement diagram. We recall that by \cite{CPS}, Th.~2.1, the necessary and sufficient
condition for a two--sided adjunction of triangulated (or abelian) categories to be a part of recollement setup is that $i_*$ induces a bijection of Hom--sets
$(X,Y)\simeq (i_*(X),i_*(Y))$. Thus the Kuhn result follows from the well known fact that $\mbox{Hom}_{{\cal P}_d}(F,G)\simeq\mbox{Hom}_{{\cal P}_{p^i d}}(F^{(i)},G^{(i)})$. On the other hand in our triangulated context we have only a monomorphism $\mbox{Ext}^*_{{\cal P}_d}(F,G)\ra\mbox{Ext}^*_{{\cal P}_{p^i d}}(F^{(i)},G^{(i)})$ (\cite{Ja},~Prop. II.10.14), hence our Derived Kan Extension  cannot be a part of recollement diagram. Nevertheless, our Theorem 3.2 computes
$i^!i_*$ in terms of something very close to the identity. This suggests an idea of enriching ${\cal DP}_d$ to turn our adjunction into a recollement.
We realize this idea in the next paper \cite{C6} where we introduce certain dg--category ${\cal P}_d^{af}$ called the category of affine strict polynomial
functors. We show there that our Derived Kan Extension ${\cal DP}_d\rec {\cal DP}_{p^i d}$ factorizes through the ``Derived Affine Kan Extension''
${\cal DP}^{af}_d\rec {\cal DP}_{p^i d}$ which is a part of recollement diagram. This explains conceptually the results of the present paper as well as
emerging some mysterious  extra structure in various Ext--computations in $\pe_d$ \cite{C2, C4, C5}.\newline

It was pointed out by Touz\'e that the Collapsing Conjecture may be used  to quickly re-obtain the Ext--computations of \cite{C2}.
We present these simplifications here in more detail than \cite{T2}; Sect. 4, Sect. 7 for the reader's convenience.
We recall from \cite{C2}, Section 5 that for
$F\in{\cal P}_d$ and a partition $\lambda$ of  $d$ we define a $d$--functor $\widetilde{F}^{\lambda}$ as the component in
$F(V_1\oplus\ldots\oplus V_d)$ of multidegree $\lambda$. Then we put $F^{\lambda}(V):=\widetilde{F}^{\lambda}(V,\ldots, V)$.
Now it follows from the Yoneda lemma (\cite{FS}, Cor. 2.12) that $\mbox{Hom}_{{\cal P}_d}(D^{\lambda}, F)\simeq F^{\lambda}({\bf k})$. Thus, by Cor.~3.7
we get \cite{C2}, Cor.~5.1
\begin{cor}
For any $F\in{\cal P}_d$ and $\lambda$ of weight $d$
\[\mbox{Ext}^*_{{\cal P}_{dp^i}}(D^{\lambda(i)},F^{(i)})\simeq F^{\lambda}(A_i).\]
\end{cor}
Also the rest of computations of \cite{C2} may be obtained with the aid of Cor.~3.7.
Namely, let $W_{\mu}, S_{\lambda}$ be respectively Weyl and
Schur functors associated to Young diagrams $\mu,\lambda$ of weight $d$, and let $s_{\mu}, s_{\lambda}$ be appropriate symmetrizations [C2, Section 3].
Then we have \cite{C2}, Th. 6.1
\begin{cor}
\[{\mbox Ext}^*_{{\cal P}_{dp^i}}(W_{\mu}^{(i)},S_{\lambda}^{(i)})\simeq s_{\mu}(s_{\lambda}(A_i^{\otimes d}\otimes{\bf k}[\Sigma_d]))\simeq
s_{\lambda}(s_{\mu}(A_i^{\otimes d}\otimes{\bf k}[\Sigma_d])).\]
\end{cor}
{\bf Proof:\ }  Cor. 3.7 allows us to replace
$\mbox{Ext}^*_{{\cal P}_{dp^i}}(W_{\mu}^{(i)},S_{\lambda}^{(i)})$ with
$\mbox{Ext}^*_{{\cal P}_d}(W_{\mu},(S_{\lambda})_{A_i}),$ i.e. in the terminology of \cite{C2} we only need to prove the ``additive version'' of the formula.
This is rather straightforward and was accomplished in the proofs of \cite{C2}; Th. 4.4, Th. 6.1.\qed
One can derive from the Collapsing Conjecture also a number of  general consequences. Among them (as pointed out by Touz\'e)
there is a simple proof of  ``twist injectivity phenomenon'' in ${\cal P}_d$. Namely, it was proved in \cite{FFSS}, Cor. 1.3 that precomposition with $I^{(i)}$ induces
a monomorphism on Ext--groups between any $F,G\in {\cal P}_d$. The proof however, went through a comparison with Ext--groups in the category of rational representations of $GL_n(\ka)$ where an analogous fact follows from  \cite{Ja}, Prop. II.10.14. Here we present a short, intrinsic proof of the twist injectivity
in ${\cal P}_d$.

\begin{cor}
For any $F,G\in {\cal P}_d$ precomposition with $I^{(i)}$ induces a monomorphism
\[{\mbox Ext}^*_{{\cal P}_d}(F,G)\longrightarrow \mbox{Ext}^*_{{\cal P}_{dp^i}}(F^{(i)},G^{(i)}).\]
\end{cor}
{\bf Proof:\ } It was observed by Touz\'e (\cite{T2}, Cor. 7.6) that the Collapsing Conjecture implies this result but, in fact, there is no need for
refering to his ``twisting spectral sequence''. We just observe that, as it is easy to see, our map fits into a commutative diagram
\[\begin{array}{ccc}
{\mbox Ext}^*_{{\cal P}_d}(F,G)&\longrightarrow& \mbox{Ext}^*_{{\cal P}_{dp^i}}(F^{(i)},G^{(i)})\\
\|&&\downarrow\\
{\mbox Ext}^*_{{\cal P}_d}(F,G)&\longrightarrow& \mbox{Ext}^*_{{\cal P}_{d}}(F,G_{A_i})
\end{array}\]
where the right vertical arrow is an isomorphism provided by Cor. 3.7 and the bottom arrow is induced by the split inclusion of functors $I\longrightarrow I\otimes A_i$.\qed
\section{Applications to nontwisted functors}
In the last section we show that the functor ${\bf K}_i$ can be effectively used in computations of $\mbox{Ext}^*_{dp^{i+j}}(F^{(i+j)},G^{(j)})$ also for nontwisted $G$.
We focus  on functors $F=W_{\mu}, G=S_{\lambda}$ for Young diagrams $\mu, \lambda$ of appropriate weights.
This class of functors
is important since  they closely approximate simple objects in ${\cal P}_{dp^i}$.
The  program of computing the groups $\mbox{Ext}^*_{{\cal P}_{dp^{i+j}}}(W_{\mu}^{(i+j)}, S_{\lambda}^{(j)})$ was started in \cite{C2} where the case
$i=0$ was handled (we re-obtained this calculation in our Cor. 3.9). In the later work \cite{C4}  some partial results  for $i>0$ were obtained.
In the present paper we  compute ${\bf K}_i(S_{\lambda})$ for certain class of $\lambda$ (Prop. 4.2) and with the aid of it the groups
$\mbox{Ext}^*_{{\cal P}_{dp^{i+j}}}(W_{\mu}^{(i+j)}, S_{\lambda}^{(j)})$ for those $\lambda$ (Cor. 4.3),  thus obtaining the Ext computations of \cite{C4} in a much simpler manner. Perhaps more importantly, we reduce in Prop. 4.1 the problem of computing $\mbox{Ext}^*_{{\cal P}_{dp^{i+j}}}(W_{\mu}^{(i+j)}, S_{\lambda}^{(j)})$ for arbitrary $\lambda$ to a very special subproblem  which should be more accessible.
\begin{prop}
Assume that
\begin{enumerate}
\item For any partition $\rho$ of  $d$, the embedding $D^{\rho(i)}\longrightarrow I^{d(i)}$ induces
an isomorphism \[\mbox{Ext}^*_{{\cal P}_{dp^i}}(D^{\rho(i)},S_{\lambda})\simeq(\mbox{Ext}^*_{{\cal P}_{dp^i}}(I^{d(i)},S_{\lambda}))_{\Sigma_{\rho}}.\]
\item ${\bf K}_i(S_{\lambda})$ is formal.
\item ${\bf K}_i(S_{\lambda})$ has  a ``good filtration'' i.e. such that  its associated object is a direct sum of Schur functors.
\end{enumerate}
Then
\[{\bf K}_i(S_{\lambda})\simeq I^{d}\otimes_{\Sigma_d}Ext^*_{{\cal P}_{dp^i}}(I^{d(i)},S_{\lambda}).\]
Moreover for any Young diagram
$\mu$ of weight $d$,
\[\mbox{Ext}^*_{{\cal P}_{dp^{i+j}}}(W_{\mu}^{(i+j)}, S_{\lambda}^{(j)})\simeq s_{\mu}(\mbox{Ext}^*_{{\cal P}_{dp^i}}(I^{d(i)},S_{\lambda})\otimes A_j^{(i)\otimes d}).\]
\end{prop}
{\bf Proof:\ }
The first part of the proposition immediately follows from Prop. 2.3.
We  now turn  to the proof of the second part.
By Cor. 3.9 and Prop. 3.1 we get
\[\mbox{Ext}^*_{{\cal P}_{dp^{i+j}}}(W_{\mu}^{(i+j)},S_{\lambda}^{(j)})\simeq
\mbox{Ext}^*_{{\cal P}_{dp^i}}(W_{\mu}^{(i)},(S_{\lambda})_{A_j})\simeq
\mbox{Ext}^*_{{\cal P}_d}(W_{\mu},{\bf K}_i((S_{\lambda})_{A_j}))\simeq\]
\[\mbox{Ext}^*_{{\cal P}_{d}}(W_{\mu},{\bf K}_i(S_{\lambda})_{A_j^{(i)}})\simeq
\mbox{Ext}^*_{{\cal P}_d}(W_{\mu},(I^{d}\otimes_{\Sigma_d}Ext^*_{{\cal P}_{dp^i}}(I^{d(i)},S_{\lambda}))_{A_j^{(i)}}).\]
Let \[0\longrightarrow D^{\mu^k}\longrightarrow\ldots\longrightarrow D^{\mu^1}\stackrel{\phi}{\longrightarrow} D^{\mu}\longrightarrow
W_{\mu}\longrightarrow 0\] be a resolution of $W_{\mu}$ by sums
of products of divided powers. To simplify notation we denote $({\bf K}_i(S_{\lambda}))_{A_j^{(i)}}\simeq(I^{d}\otimes_{\Sigma_d}Ext^*_{{\cal P}_{dp^i}}(I^{d(i)},S_{\lambda}))_{A_j^{(i)}}$
by $X$. Since $\mbox{Ext}^{>0}_{{\cal P}_d}(W_{\mu}, S_{\rho})=0$ for any $\rho$, we have $\mbox{Ext}^{>0}_{{\cal P}_d}(W_{\mu}, X)$ by
 our assumption on existence of  good filtration on ${\bf K}_i(S_{\lambda})$ and the Decomposition Formula (\cite{C2}, Cor. 2.4).
Hence, the sequence
\[0\longrightarrow \mbox{Hom}_{{\cal P}_d}(W_{\mu}, X)\longrightarrow \mbox{Hom}_{{\cal P}_d}(D^{\mu}, X)
\stackrel{\phi^*}{\longrightarrow}\ldots\longrightarrow
\mbox{Hom}_{{\cal P}_d}(D^{\mu^k}, X)\longrightarrow 0\]
is exact. Thus
\[\mbox{Hom}_{{\cal P}_d}(W_{\mu},X)=\ker(\mbox{Hom}_{{\cal P}_d}(D^{\mu}, X)
\stackrel{\phi^*}{\longrightarrow} \mbox{Hom}_{{\cal P}_d}(D^{\mu^1}, X)).\]
We describe this kernel by the arguments used in the proof of \cite{C2}, Th. 6.1.
In order to further simplify notation we put $Y:=\mbox{Hom}_{{\cal P}_d}(I^d, X)$.
Then by the fact that $X$ has  good filtration and  \cite{C2}, Lemma 6.2 we get
\[\mbox{Hom}_{{\cal P}_d}(D^{\mu}, X) \simeq (Y)_{\Sigma_{\mu}},\ \ \mbox{Hom}_{{\cal P}_d}(D^{\mu^1}, X) \simeq (Y)_{\Sigma_{\mu^1}}.\]
This, as was shown in the proof of \cite{C2}, Th. 6.1, allows us to rewrite this kernel as
\[\ker(s^{\mu}(Y)\stackrel{\phi(Y))}{\longrightarrow}s^{\mu^1}(Y))=s_{\mu}(Y)\]
(we identified here $\phi$ with the corresponding transformation between symmetrizations). Then it remains to observe that
\[Y=\mbox{Hom}_{{\cal P}_d}(I^d, (I^{d}\otimes_{\Sigma_d}Ext^*_{{\cal P}_{dp^i}}(I^{d(i)},S_{\lambda}))_{A_j^{(i)}})\simeq\]
\[\simeq\mbox{Hom}_{{\cal P}_d}(I^d,I^d\otimes_{\Sigma_d}\mbox{Ext}^*_{{\cal P}_{dp^i}}(I^{d(i)},S_{\lambda}))\otimes A_j^{(i)\otimes d}\simeq
\mbox{Ext}^*_{{\cal P}_{dp^i}}(I^{d(i)},S_{\lambda})\otimes A_j^{(i)\otimes d},\]
we use here two easy general facts:  that
$\mbox{Hom}_{{\cal P}_d}(I^d, F_U)\simeq \mbox{Hom}_{{\cal P}_d}(I^d, F)
\otimes (U)^{\otimes d}$ for any $F\in{\cal P}_d$ and $U\in{\cal V}^{gr}$, and that
$\mbox{Hom}_{{\cal P}_d}(I^d, I^d\otimes_{\Sigma_d} M)\simeq M$ for any $\Sigma_d$--module $M$. This proves the second part of the Proposition.\qed
{\bf Remark:\ }
We see that by Prop. 4.1 the problem of computing ${\bf K}_i(S_{\lambda})$ and then that of computing
$\mbox{Ext}^*_{{\cal P}_{dp^{i+j}}}(W_{\mu}^{(i+j)},S_{\lambda}^{(j)})$
is essentially reduced to  describing $\mbox{Ext}^*_{{\cal P}_{dp^i}}(I^{d(i)},S_{\lambda})$
 as a graded $\Sigma_d$-module (partial results obtained so far suggest that the hypotheses of Prop. 4.1 are satisfied for all $\lambda$). This problem for general $\lambda$ will be addressed in a future work \cite{C5}. Here we restrict our attention to certain special
case which was distinguished in \cite{C1}.\newline

Let $F_k(\lambda)$ be a Young diagram with a trivial $p$-core and the $p$-quotient with only nontrivial $k$-th diagram which is $\lambda$, and let $F^{i+1}_k(\lambda):=F_k(F^{i}_k(\lambda))$ (see \cite{C1}, Section 5).
\begin{prop}
For any Young diagram $\lambda$ of weight $d$, $i>0$, $0\leq k\leq p-1$, we have an isomorphism of graded $\Sigma_d$--modules
\[\mbox{Ext}^*_{{\cal P}_{dp^i}}(I^{d(i)},S_{F^i_k(\lambda}))\simeq Sp_{\lambda}[h^i_k],\]
where $Sp_{\lambda}$ is the Specht module associated to $\lambda$ (\cite{JK}, Chap. 7.1) and $[h^i_k]$ is a suitable shift of grading (see \cite{C4}, Th. 4.4).\newline
Moreover, $S_{F_k^i(\lambda)}$ satisfies the first assumption of Prop. 4.1, hence
\[{\bf K}_i(S_{F_k^i(\lambda)})\simeq S_{\lambda}[h_k^i].\]
\end{prop}
{\bf Remark:\ } These facts follow from \cite{C4}, Th. 4.4  but we shall prove them  independently, since, as we will see,
\cite{C4}, Th. 4.4 may be deduced from our Prop. 4.1, Prop. 4.2. This way we will obtain the Ext--computations of \cite{C4}
in a much simpler way.\newline

{\bf Proof:\ }
The proof is a rather eclectic mixture of arguments from \cite{C4} but it is still much simpler than that of \cite{C4}, Th. 4.4 (in particular we do not
refer to the Schur--de Rham complex).\newline
We need from combinatorial machinery developed in \cite{C4}  some properties of
``homological structural arrows'':
$\phi : \Lambda^{F^i_k(\lambda)}\longrightarrow S_{F^i_k(\lambda)}$, $\psi: S_{F^i_k(\lambda)}\longrightarrow S^{{\widetilde F^i_k(\lambda)}}$
(strictly speaking $\phi$ and $\psi$ exist in a suitable  localized category to which one can transport computations of the Ext--groups). It was shown in \cite{C4}, Section 3.2]  that the image of the map induced by $\psi\circ\phi$ on $\mbox{Ext}^*_{{\cal P}_{dp^i}}(I^{d(i)},-)$ equals $Sp_{\lambda}$.
 Thus to finish the proof of the first formula it suffices to observe that $\dim(\mbox{Ext}^*_{{\cal P}_{dp^i}}(I^{d(i)},S_{F^i_k(\lambda)})=\dim(Sp_{\lambda})$ (\cite{C4}, p. 46).\newline We show that for any partition $\rho$ of  $d$ the embedding $D^{\rho(i)}\longrightarrow I^{d(i)}$ induces on $\mbox{Ext}^*_{{\cal P}_{dp^i}}(-,S_{F_k^i(\lambda)})$  taking the coinvariants by induction on $d$. By the Decomposition Formula and [C4, Fact 3.4], it suffices to
 show this for $\rho=(d)$. Since  \[(\mbox{Ext}^*_{{\cal P}_{dp^i}}(I^{d(i)}, S_{F^i_k(\lambda)}))_{\Sigma_d}\simeq (Sp_{\lambda})_{\Sigma_d}[h_k^i]\simeq
 (\mbox{Hom}_{{\cal P}_{d}}(I^{d},S_{\lambda}))_{\Sigma_d}[h_k^i]\simeq\]\[\simeq\mbox{Hom}_{{\cal P}_d}(D^d,S_{\lambda})[h_k^i]\simeq S_{\lambda}({\bf k})[h_k^i],\]
by \cite{C2}, Th. 6.1, we see that $\mbox{Ext}^*_{{\cal P}_{dp^i}}(I^{d(i)}, S_{F^i_k(\lambda)}))_{\Sigma_d}\neq 0$
 if and only if $\lambda=(1^d)$. Moreover for $\lambda=(1^d)$ our formula follows from \cite{FFSS}, Th. 4.5 (strictly speaking this is the case for $k=0$, for $k>0$ the
 Ext--groups are just shifted, which may be shown by applying the Littlewood--Richardson Rule \cite{Bo} to $S_{(a,1^{dp^i-a-1})}\otimes I$). Thus it remains to show that $\mbox{Ext}^*_{{\cal P}_{dp^i}}(D^{d(i)}, S_{F^i_k(\lambda)})=0$ for $\lambda\neq (1^d)$. Let us first consider the case when $\lambda$ contains
 the diagram $(2,2)$. We will show that $\mbox{Ext}^*_{{\cal P}_{dp^i}}(D^{dp^{i-s}(s)}, S_{F^i_k(\lambda)})=0$ (and also that
$\mbox{Ext}^*_{{\cal P}_{dp^i}}(\Lambda^{dp^{i-s}(s)}, S_{F^i_k(\lambda)})=0$
 by induction on $s$. Assume these vanishings
 for $s-1$ and consider the hyperExt spectral sequences converging to $\mbox{HExt}_{{\cal P}_{dp^i}}^*({\bf dR}^{dp^{i-s+1}(s-1)},S_{F_k^i(\lambda)})$,
where ${\bf dR}^{dp^{i-s+1}(s-1)}$ is the $(s-1)$--twisted dual de Rham complex
\[0\longrightarrow \Lambda^{dp^{i-s+1}(s-1)}\longrightarrow\ldots\longrightarrow D^{dp^{i-s+1}(s-1)}\longrightarrow 0.\]
By the induction assumption and the Decomposition Formula the first spectral sequence is trivial. Hence the second spectral sequence converges to $0$.
But by the induction assumption and the Decompositon Formula all but the outer ones terms in its $E_2$--page are trivial. Thus we get the
shift in grading
\[\mbox{Ext}^*_{{\cal P}_{dp^i}}(D^{d(i)}, S_{F^i_k(\lambda)})\simeq\mbox{Ext}^*_{{\cal P}_{dp^i}}(\Lambda^{d(i)}, S_{F^i_k(\lambda)})[dp^{i-s}-1].\]
Now we  analyze in a similar manner the spectral sequence converging to
$\mbox{HExt}_{{\cal P}_{dp^i}}^*({\bf Ko}^{dp^{i-s}(s)},S_{F_k^i(\lambda)})$,
where ${\bf Ko}^{dp^{i-s}(s)}$ is the $s$--twisted dual Koszul complex
\[0\longrightarrow D^{dp^{i-s+1}(s-1)}\longrightarrow\ldots\longrightarrow \Lambda^{dp^{i-s+1}(s-1)}\longrightarrow 0.\]
From this we get the opposite shift in grading
\[\mbox{Ext}^*_{{\cal P}_{dp^i}}(\Lambda^{d(i)}, S_{F^i_k(\lambda)})\simeq
\mbox{Ext}^*_{{\cal P}_{dp^i}}(D^{d(i)}, S_{F^i_k(\lambda)})[dp^{i-s}-1],\]
which shows that the considered groups are trivial. It remains to consider the case of $\lambda$ not containing $(2,2)$ i.e. of the form
$(a,1^{d-a})$. We first take $\lambda=(2,1^{d-2})$ and  apply $\mbox{Ext}^*_{{\cal P}_{dp^i}}(D^{d(i)},-)$ to the Littlewood--Richardson
filtration \cite{Bo} on $S_{F_k^i((1^{d-1}))}\otimes S_{F_k^i((1))}$. After neglecting terms with trivial Ext--groups,  we get the long exact sequence
\[\rightarrow\mbox{Ext}^*(D^{d(i)}, S_{F^i_k((2,1^{d-2}))})\rightarrow\mbox{Ext}^*(D^{d(i)},S_{F_k^i((1^{d-1}))}\otimes S_{F_k^i((1))})\rightarrow\mbox{Ext}^*(D^{d(i)}, S_{F^i_k((1^d))})\rightarrow.\]
By \cite{FFSS}, Th. 4.5 the right arrow is an isomorphism, thus $\mbox{Ext}^*_{{\cal P}_{dp^i}}(D^{d(i)}, S_{F^i_k((2,1^{d-2})})=0$. Then by using
the Littlewood--Richardson Rule \cite{Bo} for $S_{F_k^i((1^{d-a+1}))}\otimes S_{F_k^i((a-1))}$,
we show inductively on $a$ that $\mbox{Ext}^*_{{\cal P}_{dp^i}}(D^{d(i)}, S_{F^i_k((a,1^{d-a})})=0$ for all $a\geq 2$.\newline
Thus, we have shown that
$S_{F_k^i(\lambda)}$ satisfies the first assumption of Prop. 4.1.
Hence we get
\[{\bf K}_i(S_{F^i_k(\lambda)})\simeq I^d\otimes_{\Sigma_d}Sp_{\lambda}[h^i_k]\simeq S_{\lambda}[h^i_k].\]
\qed
Now we obtain the promised Ext--computations from [C4].
\begin{cor}
For any $\mu,\lambda$ of weight $d$,  $i>0$, $0\leq k\leq p-1$,
\[\mbox{Ext}^*_{{\cal P}_{dp^{i+j}}}(W_{\mu}^{(i+j)},S_{F^i_k(\lambda)}^{(j)})\simeq s_{\mu}(s_{\lambda}(A_j^{(i)\otimes d}\otimes{\bf k}[\Sigma_d]))[h^i_k]\simeq s_{\lambda}(s_{\mu}(A_j^{(i)\otimes d}\otimes{\bf k}[\Sigma_d]))[h^i_k].\]
\end{cor}
{\bf Proof:\ }
We observe that, by Prop. 4.2,
${\bf K}_i(S_{F^i_k(\lambda)})$ has a good filtration. Thus by applying
Prop. 4.2, Prop. 4.1  we get
\[\mbox{Ext}^*_{{\cal P}_{dp^{i+j}}}(W_{\mu}^{(i+j)},S_{F^i_k(\lambda)}^{(j)})\simeq s_{\mu}(Sp_{\lambda}\otimes A_j^{(i)\otimes d}[h_k^i])\simeq
s_{\mu}(s_{\lambda}({\bf k}[\Sigma_d])\otimes A_j^{(i)\otimes d})[h_k^i]=\]
\[=s_{\mu}(s_{\lambda}({\bf k}[\Sigma_d]\otimes A_j^{(i)\otimes d}))[h_k^i].\]
The fact that $s_{\mu}$ and $s_{\lambda}$ commute in our situation is purely formal and follows e.g. from \cite{C2}, Th. 6.1.\qed
\section*{Acknowledgements}
I would like to thank Antoine Touz\'e for many discussions on functor categories, especially for detecting  a serious error in
the proof of Theorem 3.2. I am  grateful to Wilberd van der Kallen whose critical reading of the preliminary version of \cite{T3} has led to discovering the above
mentioned error and to Stanis\l aw Betley for reading  the present article.
\section*{Funding}
This work was supported by Narodowe Centrum Nauki grant 2011/01/B/ST1/06184.

\end{document}